\documentclass[12pt]{article}
\usepackage{amssymb}
\usepackage{amsfonts}
\usepackage{latexsym}
\usepackage{amsthm}

\begin{document}

%%%%%%%%%% Some definitions %%%%%%%%%%

%%%%% Equations, labels, theorems %%%%%
\renewcommand{\theequation}{\arabic{section}.\arabic{equation}}
\def\[{\begin{equation}}              \def\]{\end{equation}}
\def\lb{\label}                       \newcommand{\er}[1]{{\rm(\ref{#1})}}
\theoremstyle{plain}
\newtheorem{theorem}{\bf Theorem}[section]
\newtheorem{lemma}[theorem]{\bf Lemma}
\newtheorem{corollary}[theorem]{\bf Corollary}
\newtheorem{proposition}[theorem]{\bf Proposition}
\newtheorem{remark}[theorem]{\it Remark}
%\theoremstyle{remark}
%\newtheorem{remark}[theorem]{\bf Remark}

%%%%% Alphabet %%%%%
\def\a{\alpha}  \def\cA{{\cal A}}     \def\bA{{\bf A}}  \def\mA{{\mathscr A}}
\def\b{\beta}   \def\cB{{\cal B}}     \def\bB{{\bf B}}  \def\mB{{\mathscr B}}
\def\g{\gamma}  \def\cC{{\cal C}}     \def\bC{{\bf C}}  \def\mC{{\mathscr C}}
\def\G{\Gamma}  \def\cD{{\cal D}}     \def\bD{{\bf D}}  \def\mD{{\mathscr D}}
\def\d{\delta}  \def\cE{{\cal E}}     \def\bE{{\bf E}}  \def\mE{{\mathscr E}}
\def\D{\Delta}  \def\cF{{\cal F}}     \def\bF{{\bf F}}  \def\mF{{\mathscr F}}
\def\c{\chi}    \def\cG{{\cal G}}     \def\bG{{\bf G}}  \def\mG{{\mathscr G}}
\def\z{\zeta}   \def\cH{{\cal H}}     \def\bH{{\bf H}}  \def\mH{{\mathscr H}}
\def\e{\eta}    \def\cI{{\cal I}}     \def\bI{{\bf I}}  \def\mI{{\mathscr I}}
\def\p{\psi}    \def\cJ{{\cal J}}     \def\bJ{{\bf J}}  \def\mJ{{\mathscr J}}
\def\vT{\Theta} \def\cK{{\cal K}}     \def\bK{{\bf K}}  \def\mK{{\mathscr K}}
\def\k{\kappa}  \def\cL{{\cal L}}     \def\bL{{\bf L}}  \def\mL{{\mathscr L}}
\def\l{\lambda} \def\cM{{\cal M}}     \def\bM{{\bf M}}  \def\mM{{\mathscr M}}
\def\L{\Lambda} \def\cN{{\cal N}}     \def\bN{{\bf N}}  \def\mN{{\mathscr N}}
\def\m{\mu}     \def\cO{{\cal O}}     \def\bO{{\bf O}}  \def\mO{{\mathscr O}}
\def\n{\nu}     \def\cP{{\cal P}}     \def\bP{{\bf P}}  \def\mP{{\mathscr P}}
\def\r{\rho}    \def\cQ{{\cal Q}}     \def\bQ{{\bf Q}}  \def\mQ{{\mathscr Q}}
\def\s{\sigma}  \def\cR{{\cal R}}     \def\bR{{\bf R}}  \def\mR{{\mathscr R}}
\def\S{\Sigma}  \def\cS{{\cal S}}     \def\bS{{\bf S}}  \def\mS{{\mathscr S}}
\def\t{\tau}    \def\cT{{\cal T}}     \def\bT{{\bf T}}  \def\mT{{\mathscr T}}
\def\f{\phi}    \def\cU{{\cal U}}     \def\bU{{\bf U}}  \def\mU{{\mathscr U}}
\def\F{\Phi}    \def\cV{{\cal V}}     \def\bV{{\bf V}}  \def\mV{{\mathscr V}}
\def\P{\Psi}    \def\cW{{\cal W}}     \def\bW{{\bf W}}  \def\mW{{\mathscr W}}
\def\o{\omega}  \def\cX{{\cal X}}     \def\bX{{\bf X}}  \def\mX{{\mathscr X}}
\def\x{\xi}     \def\cY{{\cal Y}}     \def\bY{{\bf Y}}  \def\mY{{\mathscr Y}}
\def\X{\Xi}     \def\cZ{{\cal Z}}     \def\bZ{{\bf Z}}  \def\mZ{{\mathscr Z}}
\def\O{\Omega}

\def\ve{\varepsilon}   \def\vt{\vartheta}    \def\vp{\varphi}    \def\vk{\varkappa}

\def\Z{{\mathbb Z}}    \def\R{{\mathbb R}}   \def\C{{\mathbb C}}
\def\T{{\mathbb T}}    \def\N{{\mathbb N}}   \def\dD{{\mathbb D}}

%%%%% Arrows %%%%%

\def\la{\leftarrow}              \def\ra{\rightarrow}            \def\Ra{\Rightarrow}
\def\ua{\uparrow}                \def\da{\downarrow}
\def\lra{\leftrightarrow}        \def\Lra{\Leftrightarrow}

%%%%% Math signs %%%%%

\let\ge\geqslant                 \let\le\leqslant
\def\/{\over}                    \def\iy{\infty}
\def\sm{\setminus}               \def\es{\emptyset}
\def\ss{\subset}                 \def\ts{\times}
\def\pa{\partial}                \def\os{\oplus}
\def\ev{\equiv}                  \def\1{1\!\!1}
\def\iint{\int\!\!\!\int}        \def\iintt{\mathop{\int\!\!\int\!\!\dots\!\!\int}\limits}
\def\el2{\ell^{\,2}}

%%%%% Typography %%%%%

\def\lt{\biggl}                  \def\rt{\biggr}
\def\no{\noindent}               \def\ol{\overline}
\def\BBox{\hspace{1mm}\vrule height6pt width5.5pt depth0pt \hspace{6pt}}
\def\wt{\widetilde}              \def\wh{\widehat}
\newcommand{\nt}[1]{{\mathop{#1}\limits^{{}_{\,\bf{\sim}}}}\vphantom{#1}}
\newcommand{\nh}[1]{{\mathop{#1}\limits^{{}_{\,\bf{\wedge}}}}\vphantom{#1}}
\newcommand{\nc}[1]{{\mathop{#1}\limits^{{}_{\,\bf{\vee}}}}\vphantom{#1}}
\newcommand{\oo}[1]{{\mathop{#1}\limits^{\,\circ}}\vphantom{#1}}
\newcommand{\po}[1]{{\mathop{#1}\limits^{\phantom{\circ}}}\vphantom{#1}}

%%%%% Math operations %%%%%

\def\Im{\mathop{\rm Im}\nolimits}
\def\Iso{\mathop{\rm Iso}\nolimits}
\def\Ker{\mathop{\rm Ker}\nolimits}
\def\Ran{\mathop{\rm Ran}\nolimits}
\def\Re{\mathop{\rm Re}\nolimits}
\def\Tr{\mathop{\rm Tr}\nolimits}
\def\arg{\mathop{\rm arg}\nolimits}
\def\const{\mathop{\rm const}\nolimits}
\def\det{\mathop{\rm det}\nolimits}
\def\diag{\mathop{\rm diag}\nolimits}
\def\dim{\mathop{\rm dim}\nolimits}
\def\dist{\mathop{\rm dist}\nolimits}
\def\res{\mathop{\rm res}\limits}
\def\sign{\mathop{\rm sign}\nolimits}
\def\supp{\mathop{\rm supp}\nolimits}

%%%%%%%%%%% End of definitions %%%%%%%%%%

\title {The inverse Sturm-Liouville  problem with mixed
boundary conditions}

\author{Dmitri Chelkak\begin{footnote}
{Correspondence author. Dept. of Math. Analysis, Math. Mech. Faculty, St.Petersburg State
University. Universitetskij pr. 28, Staryj Petergof, 198504 St.Petersburg, Russia, e-mail:
delta4@math.spbu.ru }
\end{footnote}
and Evgeny Korotyaev\begin{footnote}
{ Institut f\"ur  Mathematik,  Humboldt Universit\"at zu
Berlin, Rudower Chaussee 25, 12489, Berlin, Germany, e-mail: evgeny@math.hu-berlin.de}
\end{footnote}
}

\maketitle

\begin{abstract}
\no Consider the operator  $H\p=-\p''+q\p=\l\p$, $\p(0)=0$, $\p'(1)+b\p(1)=0$ acting in
$L^2(0,1)$, where $q\in L^2(0,1)$ is a real potential. Let $\l_n(q,b)$, $n\ge 0$, be the
eigenvalues of $H$ and $\n_n(q,b)$ be the so-called norming constants. We give a complete
characterization of all spectral data $(\{\l_n\}_0^\iy;\{\n_n\}_0^\iy)$ that correspond to
$(q;b)\in L^2(0,1)\ts\R$. If $b$ is fixed, then we obtain a similar characterization and
parameterize the iso-spectral manifolds.
\end{abstract}

\section{Introduction and main results}

Consider  the Sturm-Liouville  problem
$$
 -\p''+q(x)\p=\l \p,\quad x\in [0,1],\ \ \ \ q\in L^2(0,1),
$$
with boundary conditions
$$
\p'(0)-a\p(0)=0,\quad \p'(1)+b \p(1)=0,\quad a,b\in \R\cup \{\iy\},
$$
where real $q$ belongs to $L^2(0,1)$, equipped with the norm $\|q\|^2=\int_0^1q(t)^2dt$.
There are a lot of papers (Borg, Gel'fand, Levitan, Marchenko, Trubowitz, ...) devoted to the
inverse spectral theory for Sturm-Liouville problems on a finite interval  (see the books
\cite{L}, \cite{M}, \cite{PT}). Recall that the inverse problem consists of the following
parts:

\no i) Prove that the spectral data (eigenvalues and some "additional"\ parameters)  {\bf
uniquely} determine the potential.

\no ii) {\bf Reconstruct} the potential from spectral data.

\no iii) {\bf Characterize} all spectral data that correspond to some fixed class of potentials.

The first uniqueness results were obtained by Borg. The first results about the reconstruction
were obtained by Gel'fand and Levitan (see details in \cite{L}, \cite{M}). Trubowitz and his
co-authors (\cite{IT}, \cite{IMT}, \cite{DT}, \cite{PT}) developed the "analytic"\ method
based on nonlinear functional analysis and the explicit reconstruction procedure for the
special case when only {\it one} spectral datum has been changed. The book \cite{PT} contains
an elegant complete solution of the inverse Dirichlet (i.e. $a\!=\!b\!=\iy$) problem on
$[0,1]$. In \cite{IT} the case when (non fixed) $a,b$ run through $\R$ was considered (we
recall the main result of \cite{IT} in Theorem \ref{TrubowitzNonFixAB}) and \cite{IMT} is
devoted to the case of {fixed} $a,b\in\R$. These authors (\cite{IT}, \cite{IMT}) reduce the
mapping $(q,a,b)\mapsto\{{\rm spectral\ data}\}$ to the case of Dirichlet boundary conditions
\cite{PT}. This reduction allows to prove the characterization theorem, to construct the
isospectral manifolds (i.e. the set of all $(q,a,b)$ with the same spectrum) and so on.

Consider the case of mixed boundary condition $a\!=\!\iy$, $b\!\in\!\R$.  We think that the
best results in this case are those of Dahlberg and Trubowitz \cite{DT}. In order to describe
their results we need some definitions. Define the self-adjoint operator $H$ in $L^2(0,1)$ by
$$
H\p=-\p''+q(x)\p,\quad \p(0)=0,\ \ \p'(1)+b\p(1)=0,\ \ b\in \R.
$$
Here and below $(\,')={\pa\/\pa x},(\dot{{\,}})={\pa\/\pa \l}$. Denote by
$\{\l_n(q,b)\}_0^\iy$ the eigenvalues of $H$.  It is well-known that all $\l_n(q,b)$ are
simple and satisfy the asymptotics
$$
\l_n(q,b)=\l_n^0+Q_0+2b+\m_n(q,b),\quad {\rm where}\quad Q_0=\int_0^1q(t)dt\ \ {\rm and}\ \
\{\m_n(q,b)\}_{0}^\iy\in\el2=\el2_0.
$$
Here and below the real Hilbert spaces $\ell_m^2$ are given by $ \ell _m^2\!=\!\left\{
\{f_n\}_0^\iy:\sum_{n\ge 0}(n\!+\!1)^{2m}f_n^2\!<\!\iy \right\}$, $m\ge 0 $, and
$\l_n^0=\pi^2(n\!+\!{1\/2})^2$, $n\ge 0$, are the unperturbed eigenvalues. The monotonicity
property $\l_0<\l_1<...$ gives that if $(q;b)$ runs through $L^2(0,1)\ts\R$, then
$\{\m_n\}_0^\iy$ doesn't run through the whole space $\el2$. In order to describe this
situation, we introduce the open and convex set
$$
\cM=\left\{\{\m_n\}_{0}^\iy\in\el2:\l_0^0\!+\!\m_0\!<\!\l_1^0\!+\!\m_1\!<\!\dots\right\}
\ss\el2.
$$
Fix some $b\in\R$. The main results of \cite{DT} are:

 \no {\it (i) The sequence of real
numbers $\l^*_n=\l_n^0+c^*+\m^*_n$, $c^*\in\R$, is the spectrum of $H$ for some potential
$q\in L^2(0,1)$ if and only if $\{\m^*_n\}_0^\iy\in\cM$.

\no (ii) Let $q\in L^2(0,1)$ and $e\in L^2_{even}(0,1)$, i.e. $e(1-x)=e(x)$, $x\in[0,1]$.
There are two cases:

\no (a) there is no potential $p\in L^2(0,1)$ such that the Dirichlet spectrum of  $p$
coincides with the Dirichlet spectrum of $e$ and $\l_n(p,b)=\l_n(q,b)$ for all $n\ge 0$; (b)
such a potential $p$ is unique.

\no An explicit condition which distinguishes the cases (a) and (b) from each other is given
in terms of $\{\l_n(q,b)\}_0^\iy$ and the Dirichlet spectrum of $e$.}

\no Roughly speaking, (ii) describes the bijection between the set of isospectral potentials
\[
\label{IsoDef} \Iso_b(\{\l_n^*\}_0^\iy)=\{q\in L^2(0,1):\l_n(q,b)=\l_n^*\ {\rm for \ all }\
n\ge 0\}.
\]
and some open subset of $L^2_{even}(0,1)$ (see \cite{DT} for details).

The main goal of our paper is to give a more explicit characterization of spectral data in
the style of \cite{PT}, \cite{IT} and to parameterize the isospectral manifolds in a more
classic way. Let $\vp(x)=\vp(x,\l,q)$ and $\x_b=\x_b(x,\l,q)$ be solutions of
$-\p''+q(x)\p=\l\p$ such that
$$
\vp(0)=0,\ \ \vp'(0)=1\quad {\rm and}\quad \x_b(1)=-1,\ \ \x'_b(1)=b.
$$
Note that the eigenvalues $\l_n(q,b)$ are the roots of the Wronskian
$$
w(\l)=w(\l,q,b)=\{\vp,\x_b\}(\l,q)=\vp'(1,\l,q)+b\vp(1,\l,q),
$$
where $\{\vp,\x_b\}=\vp\x'_b-\vp'\x_b$. The Hadamard Factorization Theorem implies the identity
\[
\label{adam}
w(\l,q,b)=\cos\sqrt\l\cdot\prod{\l-\l_n(q,b)\/\l-\l_n^0}\,,\qquad \l\in\C.
\]
Let $\p_n(x)=\p_n(x,q,b)$ be the $n$-th normalized
eigenfunction of $H$ such that $\p'_n(0)\!>\!0$. We introduce the norming constants
("additional"\ spectral data) by
\[
\label{NuDef} \n_n(q,b)=\log\left[(-1)^n\vp(1,\lambda_n(q,b),q)\right]=
\log\lt|{\p_n(1,q,b)\/\p_n'(0,q,b)}\rt|\,,\quad n\ge 0,
\]
$$
\n_n^0=\nu_n(0,0)=-\log k_n^0,\quad {\rm where}\quad
k_n^0=\sqrt{\l_n^0}=\pi(n\!+\!{\textstyle{1\/2}}).
$$
Recall that $\l_n(q,b)=\l_n^0+Q_0+2b+\m_n(q,b)$. Our main result is
\begin{theorem}\label{ThmNonFixedB} (i) The mapping
$$
\F:(q;b)\mapsto \left(Q_0\!+\!2b\,;\{\m_n(q,b)\}_{0}^\iy;\{\n_n(q,b)-\n_n^0\}_{0}^\iy\right)
$$
is a real-analytic isomorphism between $L^2(0,1)\ts\R$ and $\R\ts\cM\ts\el2_1$.\\
(ii) For each $(q;b)\in L^2(0,1)\ts\R$ the following identity is fulfilled:
\[
\label{IdentityB}  b=\sum_{n\ge 0} \lt(2-{e^{\n_n(q,b)}\/|\dot{w}(\l_n,q,b)|}\rt).
\]
\end{theorem}

\no {\it Remarks.}\ i) In the proof of (i) we use the method from \cite{PT}. The main
ingredients are nonlinear functional analysis and the explicit reconstruction procedure, when
only one $\l_n$ or $\n_n$ has been changed. Using similar arguments, it is possible to reprove
the main result of \cite{IT}, i.e. the complete characterization of spectral data in the case
$(q,a,b)\in L^2(0,1)\ts\R^2$, without the reduction to the inverse Dirichlet problem.\\
ii) Identity (\ref{IdentityB}) gives the explicit expression of $b$ in terms of spectral data
(in the case $a=\iy$). This can be rewritten in the form $b=\sum_{n\ge 0}
(2-\|\x_b(\cdot,\l_n,q)\|^{-2}_{L^2(0,1)})$ (see Lemma \nolinebreak \ref{LemmaIdentityB}).
Note that similar identities were proved in \cite{JL} using the technique of transmutation
operators and the Gel'fand-Levitan equation. Our proof is based on the contour integration. In
Sect. 5 (Appendix) we prove the analogue of (\ref{IdentityB}) in the case $a,b\in\R$.

Fix some $b\in\R$ (for instance, $b\!=\!0$ gives the boundary conditions $\p(0)\!=\!0$,
$\p'(1)\!=\!0$).  In this case spectral data $\{\l_n\}_0^\iy$, $\{\n_n\}_0^\iy$ are not
independent since they satisfy the nonlinear equation (\ref{IdentityB}). Fortunately, the
first eigenvalue $\l_0(q,b)$ can be uniquely reconstructed from the other spectral data (i.e.
$\{\l_n\}_1^\iy$ and $\{\n_n\}_0^\iy$). More precisely, we have

\begin{corollary} \label{CorFixedBEV}
For any fixed $b\in\R$ the mapping
$$
\F_b:q\mapsto \left(Q_0\,;\{\m_{n+1}(q,b)\}_{0}^\iy,\{\n_n(q,b)-\n_n^0\}_{0}^\iy\right)
$$
is a real-analytic isomorphism between $L^2(0,1)$ and $\R\ts\cM^{(1)}\ts\el2_1$, where
$$
\cM^{(1)}= \left\{\{\m_{n+1}\}_0^\iy\in\el2:
\l_1^0\!+\!\m_1\!<\!\l_2^0\!+\!\m_2\!<\!\dots\right\} \ss\el2.
$$
\end{corollary}

It is possible to "remove"\ from the spectral data not only the first eigenvalue $\l_0$ but
also one norming constant. Recall that for any sequence $\{\l_n^*\}_0^\iy$ such that
$\l_n^*=\l_n^0+c^*+\m_n^*$, where $\left(c^*;\{\m_n^*\}_0^\iy\right)\in\R\ts\cM$, the set of
isospectral potentials $\Iso_b(\{\l_n^*\}_0^\iy)$ is defined by (\ref{IsoDef}). Note that
$$
w(\l,q,b)=w^*(\l)=\cos\sqrt\l\cdot\prod_{n\ge 0}{\l\!-\!\l_n^*\/\l\!-\!\l_n^0} \quad{\rm for\
each}\ \ q\!\in\! \Iso_b(\{\l_n^*\}_0^\iy).
$$
\begin{corollary}\label{CorFixedBNC}
Fix some $b\in\R$ and let $\left(c^*;\{\m_n^*\}_0^\iy\right)\in\R\ts\cM$. Then for each $m\ge 0$
the mapping
$$
q\mapsto \{\n_n(q,b)-\n_n^0\}_{n=0,n\ne m}^{\iy}
$$
is a real-analytic isomorphism between $\Iso_b(\{\l_n^*\}_0^\iy)$ and the open nonempty set
$\cN^b_m\ss\el2_1$ given by
\[
\label{NbmDef} \cN^b_m= \lt\{\{\n_n-\n_n^0\}_{n=0,n\ne m}^\iy\in\el2_1: \sum_{n:n\ne m }
\lt(2-{e^{\n_n}\/|\dot{w}^*(\l_n^*)|}\rt)>b-2\rt\}.
\]
\end{corollary}

\no {\it Remark.}\ Let $\{\n_n-\n_n^0\}_{n:n\ne m}=t\g$, where $t\in\R$, $\g=\{\g_n\}_{n:n\ne
m} \in \el2_1$ and $\g_n> 0$, $n\ne m$. Then, due to the monotonicity reasons, there exists
$t_0\in\R$ such that $t\g\in\cN^b_m$ iff $t<t_0$.

In conclusion, note that the method from \cite{PT} works well in other inverse spectral
problems with a purely discrete spectrum. In particular, this scheme was applied in the paper
\cite{CKK2} devoted to the inverse problem for the perturbed harmonic oscillator on $\R$. The
inverse problem for the harmonic oscillator in $\R^3$ perturbed by a spherically symmetric
potential is reformulated as a problem on the half-line $\R_+$ and
it is solved in \cite{CK}. In the
last case the boundary conditions are essentially nonsymmetric and it gives an extra
motivation to investigate the case of mixed conditions on $[0,1]$, which are also essentially
nonsymmetric.

We describe the plan of the paper. In Sect. 2 we prove identity (\ref{IdentityB}) and the
Uniqueness Theorem. Sect. 3 is devoted to the analytic properties of the mapping $\F$. In
Sect. 4 we prove our main results: Theorem \ref{ThmNonFixedB} and Corollaries
\ref{CorFixedBNC}, \ref{CorFixedBEV}. Sect. 5 (Appendix) is devoted to the modifications of
identity (\ref{IdentityB}) and Corollary \ref{CorFixedBEV} to the case $a,b\in\R$.

\vskip 6pt

\no {\bf Acknowledgments.}\ Dmitry Chelkak was partly supported by grants VNP
Minobrazo\-vaniya 3.1--4733, RFFR 03--01--00377 and NSh--2266.2003.1. Evgeny Korotyaev was
partly supported by DFG project BR691/23-1. Some part of this paper was written at the
Mittag-Leffler Institute, Stockholm. The authors are grateful to the Institute for the
hospitality. The authors would like to thank Markus Klein for useful discussions.

\section{Preliminaries and proof of identity (\ref{IdentityB})}
\setcounter{equation}{0}

Introduce the fundamental solutions $\vt(x,\l,q)$, $\vp(x,\l,q)$ of the differential equation
\[
\label{DiffEq} -\psi''+q(x)\psi=\lambda\psi,\quad \l\in\C,
\]
such that
$$
\vt(0,\l,q)=\vp'(0,\l,q)=1,\qquad \vt'(0,\l,q)=\vp(0,\l,q)=1.
$$
For each $x\in [0,1]$ the functions $\vt,\vt',\vp,\vp'$ are entire with respect to
$(\l,q)\in\C\ts L_\C^2(0,1)$ (see [PT]). Moreover, the following asymptotics are fulfilled:
\[
\label{VtAsymp} \vt(x,\l,q)=\cos\sqrt{\l}x+
{1\/2\sqrt{\l}}\int_0^x\left(\sin\sqrt{\l}x+\sin\sqrt{\l}(x-2t)\right)q(t)dt+
O\lt({e^{|\Im\sqrt{\l}|x}\/|\l|}\rt),
\]
\[
\label{Vt'Asymp} \vt'(x,\l,q)=-\sqrt{\l}\sin\sqrt{\l}x+
{1\/2}\int_0^x\left(\cos\sqrt{\l}x+\cos\sqrt{\l}(x-2t)\right)q(t)dt+
O\lt({e^{|\Im\sqrt{\l}|x}\/|\l|^{1/2}}\rt),
\]
\[
\label{VpAsymp} \vp(x,\l,q)= {\sin\sqrt{\l}x\/\sqrt{\l}}+
{1\/2\l}\int_0^x\left(-\cos\sqrt{\l}x+\cos\sqrt{\l}(x-2t)\right)q(t)dt+
O\lt({e^{|\Im\sqrt{\l}|x}\/|\l|^{3/2}}\rt),
\]
\[
\label{Vp'Asymp} \vp'(x,\l,q)= \cos\sqrt{\l}x+
{1\/2\sqrt{\l}}\int_0^x\left(\sin\sqrt{\l}x-\sin\sqrt{\l}(x-2t)\right)q(t)dt+
O\lt({e^{|\Im\sqrt{\l}|x}\/|\l|}\rt)
\]
as $|\l |\to\iy$, uniformly on bounded sets of $(x;q)\in [0,1]\ts L^2_C(0,1)$.

\begin{lemma}\label{PsiN2}
For each $(q;b)\in L^2(0,1)\ts\R$ and $n\ge 0$ the following identities are fulfilled:
$$
\|\vp(\cdot,\l_n,q)\|^2=(-1)^{n+1}e^{\n_n}\dot{w}(\l_n,q,b),\qquad
\|\xi_b(\cdot,\l_n,q)\|^2=(-1)^{n+1}e^{-\n_n}\dot{w}(\l_n,q,b),
$$
$$
\p_n^2(x,q,b)={\vp(x,\l_n,q)\x_b(x,\l_n,q)\/\dot{w}(\l_n,q,b)},
$$
where $\l_n=\l_n(q,b)$ and $\n_n=\n_n(q,b)$.
\end{lemma}
\begin{proof}
Let $\vp(x)=\vp(x,\l_n,q)$, $\x_b(x)=\x_b(x,\l_n,q)$ and so on. Definition (\ref{NuDef})
yields
$$
\vp(x)=-\vp(1)\x_b(x)=(-1)^{n+1}e^{\n_n}\x_b(x).
$$
This identity and $\vp^2=\{\dot{\vp},\vp\}'$ give
$$
\|\vp\|^2=\{\dot\vp,\vp\}(1)= (-1)^{n+1}e^{\n_n}(\dot{\vp}\x'_b-\dot{\vp}'\x_b)(1)=
(-1)^{n+1}e^{\n_n}\dot{w}(\l_n,q,b),
$$
since $\dot{\vp}(0)=\dot{\vp}'(0)=\dot{\xi}_b(1)=\dot{\xi}_b'(1)=0$. Furthermore, we have
$$
\p_n^2(x)={\vp^2(x)\/\|\vp\|^2}={\vp(x)\cdot (-1)^{n+1}e^{-\n_n}\vp(x)\/\dot{w}(\l_n,q,b)}=
{\vp(x)\x_b(x)\/\dot{w}(\l_n,q,b)}
$$
and $\|\x_b\|^2=e^{-2\n_n}\|\vp\|^2=(-1)^{n+1}e^{-\n_n}\dot{w}(\l_n,q,b)$.
\end{proof}

\begin{lemma}\label{LemmaIdentityB}
For each $(q;b)\in L^2(0,1)\ts\R$ the identity
$$
b=\sum_{n\ge 0} \lt(2-{e^{\n_n}\/|\dot{w}(\l_n,q,b)|}\rt)
$$
holds true, where $\l_n=\l_n(q,b)$ and $\n_n=\n_n(q,b)$, $n\!\ge\! 0$.
\end{lemma}

\begin{proof}
Consider the meromorphic function
$$
f(\l)={\vp(1,\l,q)\/w(\l,q,b)}={\vp\/\vp'+b\vp}\,(1,\l,q),\qquad \l\in\C.
$$
All roots $\l_n=\l_n(q,b)$ of $w(\l)=w(\l,q,b)$ are simple and
$\sign\dot{w}(\l_n)=(-1)^{n+1}$. Therefore, definition (\ref{NuDef}) gives
\[
\label{xResF}
\res_{\l=\l_n}f(\l)={\vp(1,\l,q)\/\dot{w}(\l_n)}={(-1)^ne^{\n_n}\/\dot{w}(\l_n)}=
-{e^{\n_n}\/|\dot{w}(\l_n)|} \,, \qquad n\ge 0,
\]
Put $|\lambda|=\pi^2m^2\to\iy$. Then, due to (\ref{VpAsymp}), (\ref{Vp'Asymp}) and
$\int_0^1e^{ik(1-2t)}q(t)dt=o(e^{|\Im k|})$, $|k|\to\iy$, we have
$$
f(\l)= {\l^{-{1/2}}\sin\sqrt\l- \l^{-1}\cos\sqrt\l\cdot {\textstyle{1\/2}}Q_0+
o(\l^{-{1}}e^{|\Im\sqrt\l|})\/
\cos\sqrt\l+\l^{-{1/2}}\sin\sqrt\l\cdot({\textstyle{1\/2}}Q_0+b)+o(\l^{-{1/2}}e^{|\Im\sqrt\l|})}
$$
$$
={\l^{-{1/2}}\tan\sqrt\l- \l^{-1}\cdot {\textstyle{1\/2}}Q_0+ o(\l^{-1})\/
1+\l^{-{1/2}}\tan\sqrt\l\cdot({\textstyle{1\/2}}Q_0+b)+o(\l^{-{1/2}})}
$$
$$
=\l^{-{1\/2}}\tan\sqrt\l- \l^{-1}\cdot {\textstyle{1\/2}}Q_0 -\l^{-1}\tan^2\sqrt\l\cdot
({\textstyle{1\/2}}Q_0+b)+o(\l^{-1}).
$$
Let $f_0(\l)=\l^{-{1\/2}}\tan\sqrt\l$. Applying the Cauchy Theorem in the disk
$|\lambda|<\pi^2m^2$, we obtain
\[
\label{xSumRes} \sum_{n=0}^m \left(\res_{\l=\l_n}f(\l)-\res_{\l=\l_n^0}f_0(\l)\right)= -
{\textstyle{1\/2}}Q_0- ({\textstyle{1\/2}}Q_0+b)\cdot\sum_{n=0}^m \res_{\l=\l_n^0}(f_0(\l))^2.
\]
Note that
$$
f_0(\l)= {\tan\sqrt\l\/\sqrt\l}= -{2\/\l\!-\!\l_n^0} \lt(1+{\l\!-\!\l_n^0\/4\l_n^0}+
O\left((\l\!-\!\l_n^0)^2\right)\rt)\quad {\rm as}\ \ \l\to \l_n^0,\ \ n\ge 0.
$$
Hence,
$$
\res_{\l=\l_n^0}f_0(\l)= -2,\qquad \res_{\l=\l_n^0}(f_0(\l))^2=-{2\/\l_n^0},\quad n\ge 0,
$$
and
$$
\sum_{n\ge 0} \res_{\l=\l_n^0}(f_0(\l))^2= -{2\/\pi^2}\sum_{n\ge 0}{1\/(n+{1\/2})^2}=
-{2\/\pi^2}\cdot {\pi^2\/2}=-1.
$$
Substituting these identities and (\ref{xResF}) into (\ref{xSumRes}), we obtain
(\ref{IdentityB}) as $m\to\infty$.
\end{proof}

\begin{theorem}[Uniqueness Theorem]
\label{UniquenessThm} Let $\l_n(q,b)=\l_n(p,h)$ and $\n_n(q,b)=\n_n(p,h)$ for all $n\ge 0$ and
some $(q;b),(p;h)\in L^2(0,1)\ts\R$. Then $(q;b)=(p;h)$.
\end{theorem}
\begin{proof} Lemma \ref{LemmaIdentityB} immediately gives $b=h$. The rest of the proof
is standard (see \nolinebreak p. \nolinebreak 58 \nolinebreak \cite{PT}). Define functions
$$
f_1(\l,x)={F_1(\l,x)\/w(\l,q,b)},\qquad F_1(\l;x)=\vp(x,\l,p)\x'_b(x,\l,q)-
\x_b(x,\l,p)\vp'(x,\l,q),
$$
$$
f_2(\l,x)={F_2(\l,x)\/w(\l,q,b)},\qquad F_2(\l,x)=\vp(x,\l,p)\x_b(x,\l,q)-
\x_b(x,\l,p)\vp(x,\l,q).
$$
Recall that
$$
w(\l,q,b)=\vp'(1,\l,q)+b\vp(1,\l,q)=\{\vp,\x_b\}(\l,q).
$$
All roots $\l_n(q,b)=\l_n(p,b)$ of the denominator $w(\cdot,q,b)$ are simple. Moreover, all
these values are roots of the numerators $F_1$, $F_2$, since definition (\ref{NuDef}) of the
norming constant $\n_n$ yields
$$
\vp(x,\l_n,p)=(-1)^{n+1}e^{\n_n}\xi_b(x,\l_n,p),\qquad
\vp(x,\l_n,q)=(-1)^{n+1}e^{\n_n}\xi_b(x,\l_n,q).
$$
Hence, the functions $f_1$, $f_2$ are entire with respect to $\l$ for each $x\in [0,1]$. Let
$|\l|=\pi^2m^2$, $m\to\iy$. Due to asymptotics (\ref{VpAsymp}), (\ref{Vp'Asymp}), we have
$$
{w(\l,q,b)}={\cos\sqrt\l+O(|\l|^{-1/2}e^{|\Im\sqrt\l|})}=
{\cos\sqrt\l}\,\cdot(1+O(|\l|^{-1/2}).
$$
Note that
$$
\xi_b(x,\l,q)=-\vt(1-x,\l,q^{*})-b\vp(1-x,\l,q^*),
$$
where $q^*(x)=q(1-x)$, $x\in[0,1]$. Therefore, asymptotics (\ref{VtAsymp})-(\ref{Vp'Asymp})
give
$$
F_1(x,\l)= -\sin\sqrt{\l}x\cdot\sin\sqrt{\l}(1-x)+ \cos\sqrt{\l}x\cdot\cos\sqrt{\l}(1-x)+
O(|\l|^{-1/2}e^{|\Im\sqrt\l|})
$$
$$
= \cos\sqrt{\l}\cdot(1+O(|\l|^{-1/2})\qquad{\rm and}\qquad
F_2(x,\l)=O(|\l|^{-1/2}e^{|\Im\sqrt\l|}).
$$
Thus,
$$
f_1(x,\l)=1+O(|\l|^{-1/2})\ \ {\rm and}\ \ f_2(x,\l)=O(|\l|^{-1/2})\quad {\rm as}\ \
|\l|=\pi^2m^2\to\infty.
$$
Since $f_1$ and $f_2$ are entire with respect to $\l$, the maximum principle yields
$f_1(x,\l)=1$ and $f_2(x,\l)=0$ for each $x\in[0,1]$ and $\l\in\C$. In other words,
$$
(\vp(x,\l,p)-\vp(x,\l,q))\x'_b(x,\l,q)-(\x_b(x,\l,p)-\x_b(x,\l,q))\vp'(x,\l,q)=0,
$$
$$
(\vp(x,\l,p)-\vp(x,\l,q))\x_b(x,\l,q)-(\x_b(x,\l,p)-\x_b(x,\l,q))\vp(x,\l,q)=0.
$$
This gives $\vp(x,\l,p)=\vp(x,\l,q)$ for all $x\in[0,1]$ and $\l\ne \l_n(q)$, $n\ge 0$. Hence,
$p=q$\,.
\end{proof}

\section{Analyticity, asymptotics and local isomorphism}
\setcounter{equation}{0}

\begin{lemma}\label{RoucheLemma}
Let $(q;b)\in L^2_\C(0,1)\ts\C$ and $N>2(\|q\|\!+\!|b|)e^{\|q\|}$. Then $w(\l)$ has exactly
$N$ roots, counted with multiplicities, in the disc $\{\l:|\l|<\pi^2N^2\}$ and for each $n>N$
exactly one simple root $\l_n(q,b)$ in the region $\{\l:|\sqrt\l-k_n^0|<\pi\big/4\}$. There
are no other roots. Moreover, the following asymptotics is fulfilled:
\[
\label{AsympLRough} \l_n(q,b)=\l_n^0+O(1)\ \ as\ \ n\to\iy,
\]
uniformly on bounded subsets of $L^2_\C(0,1)\ts\C$.
\end{lemma}
\begin{proof}
The standard estimates of $\vp$ and $\vp'$ (see p. 13 \cite{PT}) give
$$
|\vp'(1,\l,q)-\cos\sqrt{\l}|\le |\l|^{-{1\/2}}\cdot\|q\|e^{\|q\|+|\Im\sqrt\l|},\quad
|\vp(1,\l,q)|\le |\l|^{-{1\/2}}\cdot e^{\|q\|+|\Im\sqrt\l|}.
$$
Hence,
$$
|w(\l,q,b)-\cos\sqrt{\l}\,|\le |\l|^{-{1\/2}}\cdot (\|q\|\!+\!|b|)e^{\|q\|+|\Im\sqrt\l|} <
{\textstyle{1\/2}N|\l|^{-1/2}e^{|\Im\sqrt\l|}}.
$$
Let $\l\in L_N\cup\bigcup_{n>N}l_n$, where $L_N=\{\l:|\l|\!=\!\pi^2N^2\}$ and
$l_n=\{\l:|\sqrt\l\!-\!k_n^0|\!=\!{\pi\big/4}\}$, $n>N$. Then the simple estimate
$4|\cos\sqrt\l|>e^{|\Im\sqrt\l|}$ is fulfilled (p. 27 \cite{PT}) and we have
$$
|w(\l,q,b)-\cos\sqrt{\l}\,|< 2N|\l|^{-1/2}\cdot|\cos\sqrt{\l}\,|<|\cos\sqrt{\l}\,|,\quad \l\in
L_N\cup{\textstyle\bigcup_{n>N}}l_n\,.
$$
Therefore, by Rouch\'e's Theorem, $w(\cdot,q,b)$ has as many roots as $\cos\sqrt{\l}$ in each
region bounded by these contours and the remaining unbounded domain. Furthermore, it follows
from (\ref{VpAsymp}), (\ref{Vp'Asymp}) that $0=w(\l_n,q,b)=\cos\sqrt{\l_n}+O(n^{-1})$ as
$n\to\iy$, uniformly on bounded subsets of $L^2_\C(0,1)\ts\C$. Note that the point $k_n^0$ is
the simple root of the function $\cos z$. Hence, we have $\sqrt{\l_n}=k_n^0+O(n^{-1})$ or,
equivalently, (\ref{AsympLRough}).
\end{proof}

Define the inner product in $L^2_\C(0,1)\ts\C$ by
$$
\langle(q;b),(p;h)\rangle =\int_0^1q(t)\ol{p(t)}dt+b\ol{h}.
$$
Below we write $d_{(q;b)}\s(q,b)=(f;g)$ iff the identity
$[d_{(q;b)}\s(q,b)](p;h)=\langle(f;g),(p;h)\rangle$ is fulfilled for all $(p;h)\in
L^2(0,1)\ts\R$. We use the similar notation $\pa_{(q;b)}\s(\l,q,b)$ for the partial derivative
of $\s$ with respect to $(q;b)$.
\begin{lemma} \label{GradientsLN}
(i) All functions $\l_n(q,b)$, $\m_n(q,b)$, $\n_n(q,b)$, $n\ge 0$, are real-analytic on
\linebreak $L^2(0,1)\ts\R$. Their derivatives are given by
\[
\label{LMGrad}  d_{(q(x);b)}\l_n(q,b)=\left(\p_n^2(x);\p_n^2(1)\right),\qquad
d_{(q(x);b)}\m_n(q,b)=\left(\p_n^2(x)-1;\p_n^2(1)-2\right),
\]
\[
\label{NGrad} d_{(q(x);b)}\n_n(q,b)= \left((\p_n\c_n)(x);(\p_n\c_n)(1)\right),
\]
where $\p_n(x)=\p_n(x,q,b)$ is the $n$-th normalized eigenfunction of the operator $H$ and
\[
\label{ChiDef} \c_n(x)=\c_n(x,q,b)= {\vt(x,\l_n,q)\/\p_n'(0,q,b)}-
\p_n(x,q,b)\int_0^1(\vp\vt)(t,\l_n,q)dt.
\]

\pagebreak

\no (ii) Let $(q;b)\in L^2(0,1)\ts\R$. Then, each function $\l_n$, $\m_n$, $\n_n$, has an
analytic continuation from $L^2(0,1)\ts\R$ into some complex ball $\{(p;h)\in
L^2_\C(0,1)\ts\C:\|(p-q;b-h)\|<\ve\}$, where $\ve=\ve(q,b)>0$ doesn't depend on $n$.
\end{lemma}
\begin{remark}
\label{ChiRemark} {\rm Since $\psi_n(0)=0$, we have $\psi_n(x)=\psi'_n(0)\vp(x)$. Therefore,
$\c_n(x)$ is the solution of (\ref{DiffEq}) for $\l=\l_n$ such that $\{\c_n,\p_n\}=1$ and
$\int_0^1 \p_n(x)\c_n(x)dx=0$. Note that these conditions define $\c_n$ uniquely. }
\end{remark}
\begin{proof}[Proof of Lemma \ref{GradientsLN}.]
(i) It is well-known (see p. 21 \cite{PT}) that
$$
\pa_{q(x)}\vp(1,\l,q)=\vp(x)(\vt(x)\vp(1)-\vp(x)\vt(1)),
$$
$$
\pa_{q(x)}\vp'(1,\l,q)=\vp(x)(\vt(x)\vp'(1)-\vp(x)\vt'(1)),
$$
where $\vp(x)=\vp(x,\l,q)$, $\vt(x)=\vt(x,\l,q)$ and so on. This gives
$$
\pa_{(q(x);b)}w(\l,q,b)= \left({\vphantom{|^|}}\left((\vp'\!+\!b\vp)(1)\vt(x)-
(\vt'\!+\!b\vt)(1)\vp(x)\right)\vp(x)\,;\,\vp(1)\right).
$$
Note that
$$
(\vp'\!+\!b\vp)(1)\vt(x)- (\vt'\!+\!b\vt)(1)\vp(x)= -\x_b(x),
$$
since both functions are solutions of (\ref{DiffEq}) with the same initial data at $x=1$.
Therefore,
$$
\pa_{(q(x);b)}w(\l,q,b)=\left(-\vp(x)\x_b(x)\,;\,-\vp(1)\x_b(1)\right).
$$
Recall that the function $w(\cdot,\cdot,\cdot)$ is entire and $w(\l_n(q,b),q,b)=0$. Then, the
Implicit Function Theorem and Lemma \ref{PsiN2} give
$$
d_{(q(x);b)}\l_n(q,b)=-{\pa_{(q(x);b)}w(\l_n,q,b)\/\dot{w}(\l_n,q,b)}=
\left(\p_n^2(x)\,;\,\p_n^2(1)\right).
$$
The identity $\m_n(q,b)=\l_n(q,b)-\l_n^0-\int_0^1q(t)dt-2b$ yields
$$
d_{(q(x);b)}\m_n(q,b)=d_{(q(x);b)}\l_n(q,b)-(1;2)= \left(\p_n^2(x)-1\,;\,\p_n^2(1)-2\right).
$$
Furthermore, definition (\ref{NuDef}) implies
$$
d_{(q(x);b)}\n_n(q,b)={d_{(q(x);b)}\vp(1,\l_n(q,b),q)\/\vp(1)} =
{\left(\pa_{q(x)}\vp(1,\l_n,q)\,;\,0\right)+ \dot{\vp}(1)d_{(q(x);b)}\l_n(q,b)\/\vp(1)}.
$$
Substituting $\pa_{q(x)}\vp(1,\l_n,q)$ and $d_{(q(x);b)}\l_n(q,b)$ into this formula and using
the identity $\p_n(x)=\p_n'(0)\vp(x)$, we obtain
$$
d_{(q(x);b)}\n_n(q,b)= \left(\p_n(x)\c_n(x)\,;\,\p_n(1)\c_n(1)\right),\qquad
\c_n(x)={\vt(x,\l_n,q)\/\p_n'(0,q,b)}- C_n\p_n(x,q,b),
$$
where $C_n=C_n(q,b)$ is some constant. In order to find $C_n$, note that
$\n_n(q\!+\!c,b)=\n_n(q,b)$ for all $c\in\R$. Therefore,
$$
0=\langle\pa_{q(x)}\n_n(q,b),1\rangle=
\int_0^1\p_n(x)\c_n(x)dx=\int_0^1\vp(x)\vt(x)dx-C_n\int_0^1\p_n^2(x)dx.
$$
Since $\int_0^1\p_n^2(x)dx=1$, this yields $C_n=\int_0^1\vp(x)\vt(x)dx$.\\
(ii)  The proof repeats p. 51 \cite{PT} for eigenvalues and p. 64 for norming constants.
\end{proof}

\begin{lemma}
\label{Asymptotics} (i) The following asymptotics are fulfilled:
\[
\label{AsympPsiChi} \p_n(x,q,b)=\sqrt{2}\sin k_n^0x +O(n^{-1}),\qquad \c_n(x,q,b)={\cos
k_n^0x\/ \sqrt{2}k_n^0}+O(n^{-2}), \quad n\to\infty,
\]
\[
\label{AsympLN} \m_n(q,b)=- \nh q^{(cos)}_{n+{1\/2}}+O(n^{-1}),\qquad\quad
\n_n(q,b)=\n_n^0+{\nh q^{(sin)}_{n+{1\/2}}\/2k_n^0}+O(n^{-2}),\quad n\to\infty,
\]
uniformly on bounded subsets of $[0,1]\ts L^2_\C(0,1)\ts\C$, where
$$
\nh q^{(cos)}_{n+{1\/2}}=\int_0^1\!q(x)\cos 2k_n^0xdx,\qquad \nh
q^{(sin)}_{n+{1\/2}}=\int_0^1\!q(x)\sin 2k_n^0xdx.
$$
(ii) The mapping $\F$ is real-analytic on $L^2(0,1)\ts\R$ and its Fr\'echet derivative is
given by
\[
\label{DqbPsi=} [d_{(q;b)}\F](\cdot)= \lt(\langle(1;2),\cdot\rangle; \left\{\langle
X_n,\cdot\rangle\right\}_{0}^\iy; \left\{\langle Y_n,\cdot\rangle\right\}_{0}^\iy\rt),
\]
$$
X_n=(\p_n^2(x)-1;\p_n^2(1)-2),\quad Y_n=((\p_n\c_n)(x);(\p_n\c_n)(1)),\quad n\ge 0.
$$
\end{lemma}
\begin{proof}
(i) Let $\l_n=\l_n(q,b)$. Due to (\ref{AsympLRough}), we have $\sqrt{\l_n}=k_n^0+O(n^{-1})$.
Asymptotics (\ref{VtAsymp})-(\ref{Vp'Asymp}) yields
$$
\vp(x,\l_n,q)={\sin k_n^0x\/k_n^0}+O(n^{-2}),\qquad
\|\vp(\cdot,\l_n,q)\|^2=\int_0^1\vp^2(x,\l_n,q)dx={1\/2\l_n^0}+O(n^{-3}).
$$
Hence,
$$
\p_n(x,q,b)={\vp(x,\l_n,q)\/\|\vp(\cdot,\l_n,q)\|}=\sqrt{2}\sin k_n^0x +O(n^{-1}).
$$
This gives $\p'_n(0,q,b)={\|\vp(\cdot,\l_n,q)\|^{-1}}={\sqrt{2}k_n^0+O(1)}$. Also,
$$
\vt(x,\l_n,q)=\cos k_n^0x+O(n^{-1}),\qquad \int_0^1(\vt\vp)(x,\l_n,q)dx=O(n^{-2}).
$$
Using definition (\ref{ChiDef}), we obtain the second asymptotics in (\ref{AsympPsiChi}).
Lemma \ref{GradientsLN} and (\ref{AsympPsiChi}) yield
$$
\n_n(q,b)= \int_0^1[d_{(q;b)}\m_n(tq,tb)](q;b)\,dt=
\l_n^0+\langle(2\sin^2k_n^0x\!-\!1;0),(q(x);b)\rangle+O(n^{-1}),
$$
$$
\n_n(q,b)= \n_n^0+ \int_0^1[d_{(q;b)}\n_n(tq,tb)](q;b)\,dt= \n_n^0+ {\langle(\sin k_n^0x\cos
k_n^0x ;0),(q(x);b)\rangle\/k_n^0}+O(n^{-2}).
$$
This gives (\ref{AsympLN}) since $2\sin^2k_n^0x=1-\cos 2k_n^0x$ and $\sin k_n^0x \cos
k_n^0x={1\/2}\sin 2k_n^0x$.

\no (ii) As it shown above, $\F:L^2(0,1)\ts\R\to\R\ts\cM\ts\el2_1$ is locally bounded and all
"coordinate functions"\ $\int_0^1q(t)dt\!+\!2b$, $\l_n(q,b)$, $\n_n(q,b)-\n_n^0$ are
real-analytic. Therefore (e.g., see p. \nolinebreak 138 \cite{PT}), $\F$ is a real-analytic
mapping and its Fr\'echet derivative is given by (\ref{DqbPsi=}).
\end{proof}

\begin{lemma} \label{BiOrtBasis} Let $(q;b)\in L^2(0,1)\ts\R$ and
$$
X_n=(\p_n^2(x)-1;\p_n^2(1)-2),\quad Y_n=((\p_n\c_n)(x);(\p_n\c_n)(1)),\quad n\ge 0,
$$
$$
Z_m=(-2(\p_m\c_m)'(x);(\p_m\c_m)(1)),\quad T_m=(2(\p_m^2)'(x);-\p_m^2(1)),\quad m\ge 0.
$$
Then vectors $(1;0)$, $\{Z_m\}_{0}^\iy$, $\{T_m\}_{0}^\iy$ form the biorthogonal basis to
$(1;2)$, $\{X_n\}_{0}^\iy$, $\{Y_n\}_{0}^\iy$. In other words,
$$
\begin{array}{rclrclrclc}
\langle(1;2)&\!\!\!,\!\!\!&(1;0)\rangle=1, & \langle(1;2)&\!\!\!,\!\!\!&Z_m\rangle=0, &
\langle(1;2)&\!\!\!,\!\!\!&T_m\rangle=0 & \cr \langle X_n&\!\!\!,\!\!\!&(1;0)\rangle=0, &
\langle X_n&\!\!\!,\!\!\!&Z_m\rangle=\d_{nm}, & \langle X_n&\!\!\!,\!\!\!&T_m\rangle=0,& \quad
n,m\ge 0.\cr \langle Y_n&\!\!\!,\!\!\!&(1;0)\rangle=0, & \langle
Y_n&\!\!\!,\!\!\!&Z_m\rangle=0, & \langle Y_n&\!\!\!,\!\!\!&T_m\rangle\,=\d_{nm}, &
\end{array}
$$
\end{lemma}
\begin{proof}
By definition, $<\!(1;2),(1;0)\!>\,=1$ and $<\!(1;2),Z_m\!>\,=\,<\!(1;2),T_m\!>\,=0$, since
$\p_m(0)=0$. Due to $\int_0^1\p_n^2(x)dx=1$ and $\int_0^1(\p_n\c_n)(x)dx=0$ (see Remark
\ref{ChiRemark}), we obtain $<\!X_n,(1;0)\!>\,=\,<\!Y_n,(1;0)\!>\,=0$. We consider
$<\!X_n,Z_m\!>$. The partial integration gives
$$
<\!X_n,Z_m\!>\,=-2\int_0^1(\p_n^2(x)\!-\!1)(\p_m\c_m)'(x)dx+ (\p_n^2(1)\!-\!2)(\p_m\c_m)(1)
$$
$$
=-2\int_0^1\p_n^2(x)(\p_m\c_m)'(x)dx+ (\p_n^2)(1)(\p_m\c_m)(1) =
-\int_0^1\left(\p_n^2(\p_m\c_m)'\!-\!(\p_n^2)'\p_m\c_m\right)(x)dx
$$
$$
=-\int_0^1\left(\{\p_n,\p_m\}\cdot \p_n\c_m+\p_n\p_m\cdot \{\p_n,\c_m\}\right)(x)dx.
$$
If $n\ne m$, then we have $\{\p_n,\p_m\}'=(\l_n\!-\!\l_m)\p_n\p_m$,
$\{\p_n,\c_m\}'=(\l_n\!-\!\l_m)\p_n\c_m$. Hence,
$$
<\!X_n,Z_m\!>\,=-{\{\p_n,\p_m\}\{\p_n,\c_m\}\/\l_n-\l_m}\Big|_0^1=0,\quad n\ne m,
$$
since $\{\p_n,\p_m\}(0)=\{\p_n,\p_m\}(1)=0$. If $n=m$, then $\{\p_n,\p_m\}=0$,
$\{\p_n,\c_n\}=-1$ and
$$
<\!X_n,Z_n\!>\,=\int_0^1\p_n^2(x)dx=1,\quad n\ge 0.
$$
The proof of other identities is similar.
\end{proof}

\begin{theorem}[Local Isomorphism]
\label{ThmLocIso} \label{DqbPsiInvert} For each $(q;b)\in L^2(0,1)\ts\R$ operator
$d_{(q;b)}\F: L^2(0,1)\ts\R\to\R\ts\el2\ts\el2_1$ given by (\ref{DqbPsi=}) is invertible.
\end{theorem}
\begin{proof} Due to Lemma \ref{Asymptotics}, we have
$$
X_n=(-\cos\pi(2n\!+\!1)x;0)+O(n^{-1}),\quad 2k_n^0Y_n=(\sin\pi(2n\!+\!1)x;0)+O(n^{-1}),\quad
n\ge 0.
$$
Note that  vectors $(0;1)$, $\{(2\cos\pi(2n\!+\!1)x;0)\}_{0}^\iy$,
$\{(2\sin\pi(2n\!+\!1)x;0)\}_{0}^\iy$ form the orthonormal basis in $L^2(0,1)\ts\R$ and the
error terms are square summable. Then, $d_{(q(x);b)}\F$ is a Fredholm operator. Due to Lemma
\nolinebreak \ref{BiOrtBasis}, the vectors $(1;2)$, $\{X_n\}_0^\iy$, $\{Y_n\}_0^\iy$ are
linearly independent. Using the standard arguments from the functional analysis (e.g., see
p.163 \cite{PT}), we deduce that $(d_{(q;b)}\F)^{-1}$ is bounded.
\end{proof}

\section{Proof of Theorem \ref{ThmNonFixedB} and
Corollaries \ref{CorFixedBNC}, \ref{CorFixedBEV}}
\setcounter{equation}{0}

\begin{lemma}\label{SurjectionSpectrum}
For each $(c^*,\{\m_n^*\}_{0}^\iy)\!\in\!\R\ts\cM$ there exists a potential $q^*\!\in\!
L^2(0,1)$ such that $\l_n(q^*,0)=\l_n^0\!+\!c^*\!+\!\m_n^*$ for all $n\!\ge\! 0$.
\end{lemma}
\begin{proof}
Define $\l_{n+{1\/2}}^0$ and $\l_{n+{1\/2}}^*$ by
$$
\l_{n+{1\/2}}^0=\pi^2(n\!+\!1)^2,\qquad \l_{n+{1\/2}}^*=
{\l_{n+1}^0\!-\!\l_{n+{1\/2}}^0\/\l_{n+1}^0\!-\!\l_{n}^0}\cdot\l_n^*+
{\l_{n+{1\/2}}^0\!-\!\l_{n}^0\/\l_{n+1}^0\!-\!\l_{n}^0}\cdot\l_{n+1}^*,\quad n\ge 0.
$$
Then, $\l_n^*<\l_{n+{1\/2}}^*<\l_{n+1}^*$, $n\ge 0$, and
$\l_{m/2}^*=\pi^2(m\big/2)^2+c^*+\m_{m/2}^*$, where $\{\m_{m/2}^*\}_0^\iy\in\el2$. Using
\cite{PT}, we see that there exists an "even"\ potential $q\in L^2(0,2)$,
$q(x)\!=\!q(2\!-\!x)$, $x\!\in\![0,2]$, such that $\{\l_{m/2}^*\}_{0}^\iy$ is the Dirichlet
spectrum of $q$ on the interval $[0,2]$. Put $q^*=q|_{[0,1]}$. Then we have
$\l_n(q^*,0)=\l_{2n/2}^*=\l_n^*$ for all $n\ge 0$.
\end{proof}

\begin{lemma}\label{DarbouxNC} Fix any $(q;b)\in L^2(0,1)\ts\R$, $n\ge 0$ and $t\in\R$. Denote
$$
q_n^t(x)=q(x)-2{d^2\/dx^2}\log\e_n^t(x,q,b), \ \ {where} \ \ \e_n^t(x,q,b)=
1+(e^t-1)\int_x^1\p_n^2(t,q,b)dt,\ \ x\in [0,1],
$$
$$
b_n^t=b-(e^t-1)\p_n^2(1,q,b).
$$
Then $(q_n^t;b_n^t)\in L^2(0,1)\ts\R$ and
$$
\l_m(q_n^t,b_n^t)=\l_m(q,b),\quad\n_m(q_n^t,b_n^t)=\n_m(q,b)+t\d_{mn}\quad {for\ all\ }m\ge 0.
$$
\end{lemma}
\begin{proof} Let $\p_n(x)=\p_n(x,q,b)$ and so on.
Repeating the arguments of p. 91--93 \cite{PT} or using direct calculations, it is easy to
check that for each $m\ge 0$ the function
$$
\wt{\p}_m(x)=\p_m(x)-(e^t-1){\p_n(x)\/\e_n^t(x)}\int_x^1\p_n(t)\p_m(t)dt,
$$
is some solutions of $-\p''+q_n^t(x)\p=\l_m(q,b)\p$ (in particular,
$\wt{\p}_n(x)={\p_n(x)\big/\e_n^t(x)}$). Since $\p_m(0)=\p_n(0)=0$, we have $\wt\p_m(0)=0$,
$m\ge 0$. The crucial point is the new boundary condition at $x=1$. We have
$$
\wt{\p}'_m(1)+b_n^t\wt{\p}_m(1)= \p'_m(1)+(e^t-1)\p_n^2(1)\p_m(1)+b_n^t\wt{\p}_m(1) =
\p'_m(1)+b\p_m(1)=0
$$
for all $m\ge 0$. Therefore, $\l_m(q_n^t,b_n^t)=\l_m(q,b)$, $m\ge 0$ (there are no other roots
due to Lemma \ref{RoucheLemma}). Furthermore,
$$
\wt{\p}_m'(0)=\p_m'(0),\ \ m\ne n,\ \ {\rm and}\ \ \wt{\p}_n'(0)=e^{-t}\wt{\p}_n'(0).
$$
By definition (\ref{NuDef}), it gives $\n_m(q_n^t,b_n^t)=\n_m(q,b)+t\d_{mn}$ for all $m\ge 0$.
\end{proof}

\begin{proof}[{\bf Proof of Theorem \ref{ThmNonFixedB}}]
ii) Identity (\ref{IdentityB}) was proved in Lemma \ref{LemmaIdentityB}.

\no (i)  It follows from Lemma \ref{Asymptotics} that $\F$ maps  $L^2(0,1)\ts\R$ into $
\R\ts\cM\ts\el2_1$ and $\F$ is real-analytic and its Fr\'echet derivative is given by
(\ref{DqbPsi=}). Theorem \ref{DqbPsiInvert} yields that $\F$ is a local real-analytic
isomorphism. Furthermore, due to Theorem \ref{UniquenessThm}, $\F$ is one-to-one. We prove
that $\F$ is onto.

Let $c^*\in\R$, $\{\m_n^*\}_0^\iy\in\cM$ and $\{\n_n^*-\n_n^0\}_0^\iy\in\el2_1$ be an
arbitrary spectral data. Due to Lemma \ref{SurjectionSpectrum}, there exists a potential
$q^*\in L^2(0,1)$ such that $\l_n(q^*,0)=\l_n^*=\l_n^0+c^*+\m_n^*$ for all $n\!\ge\!0$. Let
$\e_n=\n_n(q^*,0)-\n_n^0$. Note that we have the following convergence:
$$
(\e_0,...,\e_{N-1},\n_N^*-\n_N^0,\n_{N+1}^*-\n_{N+1}^0,...)\to (\e_0,...,\e_{N-1}, \e_N,...)\
\ {\rm in}\ \ \el2_1\ \ {\rm as}\ \ N\to\iy.
$$
Theorem \ref{ThmLocIso} yields that the mapping $\F$ is invertible in some neighborhood of the
point $\F(q^*,0)=(c^*;\{\m_n^*\}_0^\iy;\{\e_n\}_0^\iy)\in\R\ts\cM\ts\el2_1$. Hence, for some
large integer $N$ and $(q^{(N)};b^{(N)})\in L^2(0,1)\ts\R$ we have
$$
\lt(c^*;\{\m_n^*\}_{0}^\iy;(\e_0,...,\e_{N-1},
\n_N^*-\n_N^0,\n_{N+1}^*-\n_{N+1}^0,...)\rt)=\F(q^{(N)},b^{(N)}).
$$
Applying Lemma \ref{DarbouxNC} step by step, we construct the sequence of potentials and
boundary constants
$$
(q^{(k)};b^{(k)})=(q^{(k+1)};b^{(k+1)})_{k}^{t_k}\in L^2(0,1)\ts\R,\ \
t_k=\n_k^*\!-\!\n_k^0\!-\!\e_k,\quad k=N\!-\!1,N\!-\!2,...,0,
$$
such that $\F(q^{(k)},b^{(k)})= \left(c^*;\{\m_n^*\}_{0}^\iy;(\e_0,...,\e_{k-1},
\n_k^*-\n_k^0,\n_{k+1}^*-\n_{k+1}^0,...)\right)$. In particular,
$\l_n(q^{(0)},b^{(0)})=\l_n^*$ and $\n_n(q^{(0)},b^{(0)})=\n_n^*$ for all $n\ge 0$. We are
done.
\end{proof}

\begin{proof}[{\bf Proof of Corollary \ref{CorFixedBNC}}]
Fix some $b,c^*\in\R$, $\{\m_n^*\}_0^\iy\in\cM$ and $m\ge 0$. It follows from Theorem
\nolinebreak \ref{ThmNonFixedB} that the mapping $q\mapsto \{\n_n(q,b)-\n_n^0\}_{0}^\iy$ is a
real-analytic isomorphism between $\Iso_b(\{\l_n^*\}_0^\iy)$ and the set
$$
\cN^b=\lt\{\{\n_n-\n_n^0\}_0^\iy\in\el2_1: \sum_{n\ge 0}
\lt(2-{e^{\n_n}\/|\dot{w}^*(\l_n^*)|}\rt) = b\rt\}\ss\el2_1.
$$
Introduce the mapping $P_{m}:\{\n_n-\n_n^0\}_0^\iy\mapsto \{\n_n-\n_n^0\}_{n=0,n\ne m}^\iy$.
Since $e^{\n_m}>0$, we have
$$
P_m: \cN^b \to \cN^b_m,
$$
where $\cN^b_m$ is given by (\ref{NbmDef}). Moreover, this mapping is a bijection between
$\cN^b$ and $\cN^b_m$ since $\n_m$ is uniquely reconstructed from $\{\n_n\}_{n:n\ne m}$ by
$$
\n_m=\log|\dot{w}^*(\l^*_0)|+\log\lt[2-b+\sum_{n:n\ne m}
\lt(2-{e^{\n_n}\/|\dot{w}(\l^*_n)|}\rt)\rt].
$$
Note that Theorem \ref{ThmNonFixedB} yields that the mapping
$$
\{\n_n(q,b)-\n_n^0\}_{0}^\iy\mapsto
\left(c^*;\{\m_n^*\}_0^\iy;\{\n_n(q,b)-\n_n^0\}_{0}^\iy\right)=\F(q,b)\mapsto b=\sum_{n\ge 0}
\lt(2-{e^{\n_n}\/|\dot{w}^*(\l_n^*)|}\rt)
$$
is real-analytic. It is clear that $\pa b\big/\pa\n_m=-|\dot{w}(\l_m)|^{-1}\ne 0$. Using the
Implicit Function Theorem, we obtain that $P_m^{-1}:\cN^b_m\to\cN^b$ is real-analytic.
Therefore, $P_m$ is a real-analytic isomorphism and the mapping $q\mapsto
\{\n_n(q,b)\!-\!\n_n^0\}_{n=0,n\ne m}^\iy$ is a composition of real-analytic isomorphisms.
\end{proof}

For each $c\in\R$ and $\{\m_n\}_0^\iy\in\cM$, we introduce the function
\[
\label{xWDef} W(\l)=W(\l\,;\,c,\m_0,\m_1,...)=\cos\sqrt\l\cdot\prod_{k\ge 0}
{\l\!-\!\l_k\/\l\!-\!\l_k^0}, \quad {\rm where}\quad \l_n=\l_n^0+c+\m_n,\ n\ge 0.
\]
\begin{lemma}
\label{LmMonot} Let $c\in\R$, $\{\m_n\}_0^\iy\in\cM$ and $n\ge 0$. Then,
$$
{\pa\/\pa\m_0}\,|\dot{W}(\l_n)|^{-1}>0,\qquad \lim_{\m_0\to-\iy}|\dot{W}(\l_n)|^{-1}=0.
$$
Moreover, $\lim_{\m_0\to\l_1^0-\l_0^0+\m_1}|\dot{W}(\l_0)|^{-1}=
\lim_{\m_0\to\l_1^0-\l_0^0+\m_1}|\dot{W}(\l_1)|^{-1}=+\iy$.
\end{lemma}
\no {\it Remark.\ } Note that the convergence $\m_0\to\l_1^0-\l_0^0+\m_1$ is equivalent to
$\l_0\to\l_1$.
\begin{proof} In order to consider $W$ as a function of $\m_0$, we introduce the notation
$$
W_\b(\l)=W(\l\,;\,c,\m_0\!+\!\b,\m_1,\m_2,...)={\l\!-\!\l_0\!-\!\b\/\l\!-\!\l_0}\cdot W(\l)=
\lt(1-{\b\/\l\!-\!\l_0}\rt)W(\l).
$$
Let $n\ge 1$. Due to $W(\l_n)=0$ and $\dot{W}(\l_n)\ne 0$, we have
$$
\dot{W}_\b(\l_n)=\lt(1-{\b\/\l_n\!-\!\l_0}\rt)\dot{W}(\l_n),\quad\quad
{\pa\dot{W}(\l_n)\/\pa\m_0}= {\pa\/\pa\b}\,\dot{W_\b}(\l_n)\Big|_{\b=0}=
-{\dot{W}(\l_n)\/\l_n\!-\!\l_0}\,.
$$
Hence,
$$
{\pa\/\pa\m_0}\,|\dot{W}(\l_n)|^{-1}= |\dot{W}(\l_n)|^{-2}\cdot
{|\dot{W}(\l_n)|\/\l_n\!-\!\l_0}>0,
$$
$$
\lim_{\m_0\to-\iy}|\dot{W}(\l_n)|^{-1}= \lim_{\b\to-\iy}|\dot{W}_\b(\l_n)|^{-1}=
\lim_{\b\to-\iy}\lt|1-{\b\/\l_n\!-\!\l_0}\rt|^{-1}\cdot|\dot{W}(\l_n)|^{-1}=0.
$$
Moreover, if $n=1$, then
$$
\lim_{\m_0\to \l_1^0-\l_0^0+\m_1}|\dot{W}(\l_1)|^{-1}= \lim_{\b\to \l_1-\l_0}|
\dot{W}_\b(\l_1)|^{-1}= \lim_{\b\to \l_1-\l_0}\,
\lt|1-{\b\/\l_1\!-\!\l_0}\rt|^{-1}\!\!\cdot|\dot{W}(\l_1)|^{-1}=+\iy.
$$
Let $n=0$. We have
$$
\dot{W}_\b(\l_0\!+\!\b)={W(\l_0\!+\!\b)\/\b}\,,\quad\quad {\pa\dot{W}(\l_0)\/\pa\m_0}=
{\pa\/\pa\b}\,\dot{W}_\b(\l_0\!+\!\b)\Big|_{\b=0}= {\ddot{W}(\l_0)\/2}\,.
$$
Recall that $\dot{W}(\l_0)<0$. Therefore, the Hadamard factorization yields
$$
W(\l)=C(\l\!-\!\l_0)\prod_{n\ge 1}\lt(1-{\l-\l_0\/\l_n\!-\!\l_0}\rt),\quad C<0,
$$
and $\ddot{W}(\l_0)=-2C\sum_{n\ge 1}(\l_n\!-\!\l_0)^{-1}>0$. Hence,
$$
{\pa\/\pa\m_0}\,|\dot{W}(\l_0)|^{-1}= -{\pa\/\pa\m_0}\,(\dot{W}(\l_0))^{-1}=
{\ddot{W}(\l_0)\/2(\dot{W}(\l_0))^2}>0,
$$
$$
\lim_{\m_0\to-\iy}|\dot{W}(\l_0)|^{-1}= \lim_{\b\to-\iy}|{W}_\b(\l_0\!+\!\b)|^{-1}=
\lim_{\b\to-\iy}|\b|\cdot|{W}(\l_0\!+\!\b)|^{-1}=0,
$$
and
$$
\lim_{\m_0\to \l_1^0-\l_0^0+\m_1}|\dot{W}(\l_0)|^{-1}=
\lim_{\b\to\l_1-\l_0}|{W}_\b(\l_0\!+\!\b)|^{-1}= \lim_{\b\to \l_1-\l_0}\,
|{W}(\l_0\!+\!\b)|^{-1}\cdot|\b|=+\iy,
$$
since $|{W}(\l_1)|=0$.
\end{proof}

For each $c\in\R$, $\{\m_n\}_0^\iy\in\cM$ and $\{\n_n-\n_n^0\}_{0}^\iy\in\el2_1$, we put
\[
\label{xBDef} B(c,\m_0,\n_0,...)= \sum_{n\ge 0} \lt(2-{e^{\n_n}\/|\dot{W}(\l_n)|}\rt),
\]
where the function $W$ is given by (\ref{xWDef}). Note that the sum converges due to Theorem
\ref{ThmNonFixedB}.

\begin{proof}[{\bf Proof of Corollary \ref{CorFixedBEV}}]
Fix some $b\in\R$ and $m\ge 0$. Due to Theorem \ref{ThmNonFixedB}, the mapping
$$
q\mapsto\F(q,b)=\left(Q_0\!+\!2b\,;\{\m_n(q,b)\}_{0}^\iy;\{\n_n(q,b)-\n_n^0\}_{0}^\iy\right)
$$
is a bijection between $L^2(0,1)$ and the set
$$
\cS^b=\lt\{(c;\{\m_n\}_0^\iy;\{\n_n-\n_n^0\}_0^\iy)\in\R\ts\cM\ts\el2_1:
B(c,\m_0,\n_0,...)\!=\!b\rt\}.
$$
Define the mapping $P:\cS^b\ss\R\ts\cM\ts\el2_1\to\R\ts\cM^{(1)}\ts\el2_1$ by
$$
(c;\{\m_n\}_0^\iy;\{\n_n\!-\!\n_n^0\}_0^\iy)\mapsto
(c;\{\m_{n+1}\}_0^\iy;\{\n_n\!-\!\n_n^0\}_{0}^\iy).
$$
Note that $\F_b(q)=P(\F(q,b))$. In particular, $\F_b$ is real-analytic as a composition of
real-analytic mappings. Due to (\ref{xBDef}) and Lemma \ref{LmMonot}, we have
$$
{\pa B(c,\m_0,\n_0,...)\/\pa\m_0}<0,\qquad \lim_{\m_0\to-\iy} B(c,\m_0,\n_0,...)=+\iy,\quad
\lim_{\m_0\to\l_1^0-\l_0^0+\m_1} B(c,\m_0,\n_0,...)=-\iy.
$$
Therefore, for each
$(c;\{\m_{n+1}\}_0^\iy;\{\n_n\!-\!\n_n^0\}_{0}^\iy)\in\R\ts\cM^{(1)}\ts\el2_1$ there exists
unique point $\m_0$ such that $\m_0<\l_1^0-\l_0^0+\m_1$ and $B(c,\m_0,\n_0,...)=b$. Thus, the
mapping $P$ is a bijection between $\cS^b$ and $\R\ts\cM^{(1)}\ts\el2_1$. Hence, $\F_b$ is a
bijection between $L^2(0,1)$ and $\R\ts\cM^{(1)}\ts\el2_1$. Using Theorem \ref{ThmNonFixedB},
we see that the mapping
$$
(c,\m_0,\n_0,...)\mapsto B(c,\m_0,\n_0,...)
$$
is real-analytic and $\pa B\big/\pa\m_0<0$ due to Lemma \ref{LmMonot}. Then, the Implicit
Function Theorem yields that $P^{-1}:\R\ts\cM^{(1)}\ts\el2_1\to\cS^b$ is real-analytic too. We
deduce that $(\F^b)^{-1}$ is real-analytic as a composition of real-analytic mappings.
\end{proof}

\section{Appendix. The case of general boundary conditions}
\setcounter{equation}{0}

Consider  the Sturm-Liouville  problem
$$
 -\p''+q(x)\p=\l \p,\quad x\in [0,1],\quad q\in L^2(0,1),
$$
$$
\p'(0)-a\p(0)=0,\quad \p'(1)+b \p(1)=0,\quad a,b\in \R.
$$
Denote by $\s_n(q,a,b)$ the eigenvalues of this problem. It is well-known that
\[
\label{SnDef} \s_n(q,a,b)=\s_n^0+Q_0+2a+2b+\t_n(q,a,b),\quad {\rm where} \quad
Q_0=\int_0^1q(t)dt,\quad \{\t_n\}_0^\iy\in\el2
\]
and $\s_n^0=\pi^2n^2$, $n\ge 0$. Note the $\s_n(q,a,b)$ are the roots of the Wronskian
$$
w(\l,q,a,b)= (\vt'\!+\!a\vp'+b(\vt\!+\!a\vp))(1,\l,q).
$$
Following \cite{IT}, we introduce the norming constants
$$
\vk_n(q,a,b)=\log[(-1)^n(\vt\!+\!a\vp)(1,\s_n,q)]=
\log\left|\p_n(1,q,a,b)\/\p_n(0,q,a,b)\right|,
$$
where $\p_n$ is the $n$-th normalized eigenfunction such that $\p_n(0)>0$.
\begin{theorem}[Isaacson, Trubowitz \cite{IT}]
\label{TrubowitzNonFixAB} The mapping
$$
\P:(q;a;b)\mapsto
\left(Q_0\!+\!2a\!+\!2b\,;\{\t_n(q,a,b)\}_{0}^\iy;\{\vk_n(q,a,b)\}_{0}^\iy\right)
$$
is a real-analytic isomorphism between $L^2(0,1)\ts\R^2$ and $\R\ts\cT\ts\el2_1$, where
$$
\cT=\left\{\{\t_n\}_{0}^\iy\in\el2:\s^0_0+\t_0\!<\!\s_1^0+\t_1\!<\!\dots\right\}\ss\el2.
$$
\end{theorem}
\no Theorem \ref{TrubowitzNonFixAB} and Theorem \ref{ThmNonFixedB} are similar. In order to
consider the case of fixed $a,b\in\R$, we need some modification of identity
(\ref{IdentityB}).
\begin{proposition}
\label{PropIdentitiesAB} For each $(q;a;b)\in L^2(0,1)\ts\R^2$ the following identities are
fulfilled:
\[
\label{IdentitiesAB} -1+\sum_{n\ge 0} \left(2-{e^{\vk_n}\/|\dot{w}(\s_n)|}\right)=b; \quad
-1+\sum_{n\ge 0} \left(2-{e^{-\vk_n}\/|\dot{w}(\s_n)|}\right)=a,
\]
where $\displaystyle w(\l)\!=\!w(\l,q,a,b)\!=\!-\sqrt\l\sin\sqrt\l\,\cdot {\prod}_{n\ge 0}
{{\l\!-\!\s_n\/\l\!-\!\s_n^0}}$ and $\s_n=\s_n(q,a,b)$, $\vk_n=\vk_n(q,a,b)$.
\end{proposition}

\no {\it Remark.}\ Note that identities (\ref{IdentitiesAB}) and Theorem
\ref{TrubowitzNonFixAB} directly yield the following well-known result:

\begin{quotation}
{\no \it Fix some constant $a=b$ and $q\in L^2_{even}(0,1)$, i.e. $q(x)=q(1-x)$, $x\in [0,1]$.
\\ If $p\in L^2(0,1)$ is such that $\s_n(p,a,a)=\s_n(q,a,a)$ for all $n\ge 0$, then
$p=q$.}
\end{quotation}

\no Indeed, in this case we have $w(\l,p,a,a)\!=\!w(\l,q,a,a)$ and $\vk_n(q,a,a) \!= \!0$ for
all $n\!\ge\! 0$. Summing identities (\ref{IdentitiesAB}) for the potential $p$ and
subtracting corresponding ones for the potential $q$, we deduce that
$$
\sum_{n\ge 0}{\cosh \vk_n(p,a,a)-1\/|\dot{w}(\s_n)|}=0.
$$
This yields $\vk_n(p,a,a)=0$ for all $n\ge 0$ since $\cosh\vk_n\ge 1$. Thus, Theorem
\ref{TrubowitzNonFixAB} gives $p=q$.

\pagebreak

Let $a,b\in\R$ be fixed. Then, due to identities (\ref{IdentitiesAB}), it is possible to
reconstruct $(\s_0,\vk_0)$ from the other spectral data $\{\s_n\}_1^\iy$, $\{\vk_n\}_1^\iy$.
More precisely, we have
\begin{proposition}
\label{CorFixedAB} For any fixed $a,b\in\R$ the mapping
$$
\P_{a,b}:q\mapsto
\left(Q_0\!+\!2a\!+\!2b\,;\{\t_{n+1}(q,a,b)\}_{0}^\iy;\{\vk_{n+1}(q,a,b)\}_{0}^\iy\right)
$$
is a real-analytic isomorphism between $L^2(0,1)$ and $\cT^{(1)}\ts\el2_1$, where
$$
\cT^{(1)}= \left\{\{\t_{n+1}\}_0^\iy\in\el2:
\s_1^0\!+\!\t_1\!<\!\s_2^0\!+\!\t_2\!<\!\dots\right\} \ss\el2.
$$
\end{proposition}

\begin{proof}[{\bf Proof of Proposition \ref{PropIdentitiesAB}.}]
We prove the first identity in (\ref{IdentitiesAB}), the proof of the second is similar.
Consider the meromorphic function
$$
f(\l)={(\vt+a\vp)(1,\l,q)\/w(\l,q,a,b)}=
{\vt+a\vp\/(\vt'\!+\!a\vp')+b(\vt\!+\!a\vp)}\,(1,\l,q),\qquad \l\in\C.
$$
All roots $\s_n=\s_n(q,a,b)$ of $w(\l)$ are simple. Therefore, definition (\ref{SnDef}) yields
\[
\label{xResFNN} \res_{\l=\s_n}f(\l)={(\vt+a\vp)(1,\l,q)\/\dot{w}(\s_n)}=
{(-1)^ne^{\vk_n}\/\dot{w}(\s_n)}=-{e^{\vk_n}\/|\dot{w}(\s_n)|}\,.
\]
Put $|\lambda|=\pi^2(m\!+\!{1\/2})^2\to\iy$. Then, due to asymptotics
(\ref{VtAsymp})-(\ref{Vp'Asymp}), we have
$$
f(\l)= {\cos\sqrt\l+\l^{-{1/2}}\sin\sqrt\l\cdot({\textstyle{1\/2}}Q_0+a)+
o(\l^{-{1/2}}e^{|\Im\sqrt\l|})\/
-\l^{1/2}\sin\sqrt\l+\cos\sqrt\l\cdot({\textstyle{1\/2}}Q_0+a+b)+o(e^{|\Im\sqrt\l|})}
$$
$$
=-\l^{-{1\/2}}\cot\l-\l^{-1}\cdot({\textstyle{1\/2}}Q_0+a)-
\l^{-1}\cot^2\l\cdot({\textstyle{1\/2}}Q_0+a+b)+o(\l^{-1}).
$$

\no Let $f_0(\l)=-\l^{-{1\/2}}\cot\sqrt\l$. Applying the Cauchy Theorem in the disk
$|\lambda|<\pi^2(m\!+\!{1\/2})^2$, we obtain
\[
\label{xSumResNN} \sum_{n=0}^m \left(\res_{\l=\s_n}f(\l)-\res_{\l=\s_n^0}f_0(\l)\right)=
-({\textstyle{1\/2}}Q_0+a)- ({\textstyle{1\/2}}Q_0+a+b)\cdot\sum_{n=0}^m
\res_{\l=\s_n^0}(f_0(\l))^2.
\]
Note that
$$
f_0(\l)= -{1\/\sqrt\l\tan\sqrt\l}= -{1\/\l}\lt(1-{\l\/3}+O(\l^2)\rt),\quad \l\to \s_0\!=\!0,
$$
and
$$
f_0(\l)=-{2\/\l\!-\!\s_n^0}\lt(1-{\l\!-\!\s_n^0\/4\s_n^0}+O((\l\!-\!\s_n^0)^2)\rt),\quad
\l\to\s_n^0,\ \ n\ge 1.
$$
Hence,
$$
\res_{\l=\s_0^0}f_0(\l)=-1,\quad \res_{\l=\s_n^0}f_0(\l)=-2,\ \ n\ge 1,
$$
$$
\res_{\l=\s_0^0}(f_0(\l))^2=-{2\/3},\quad \res_{\l=\s_n^0}(f_0(\l))^2=-{2\/\s_n^0},\quad n\ge
1,
$$
and
$$
\sum_{n\ge 0} \res_{\l=\s_n^0}(f_0(\l))^2= -{2\/3}-{2\/\pi^2}\sum_{n\ge 0}{1\/n^2}=
-{2\/3}-{2\/\pi^2}\cdot {\pi^2\/6}=-1.
$$
Substituting these identities and (\ref{xResFNN}) into (\ref{xSumResNN}), we obtain
(\ref{IdentitiesAB}) as $m\to\iy$.
\end{proof}

For each $c\in\R$, $\{\t_n\}_{0}^\iy\in\cT$ and $\{\vk_n\}_{0}^\iy\in\el2_1$ we put
$$
W(\l)=W(\l;c,\t_0,\t_1,...)= -\sqrt\l\cos\sqrt\l\cdot \prod_{k\ge 0}
{\l-\s_k\/\l-\s_k^0},\quad {\rm where}\quad \s_k=\s_k^0+c+\t_k,\ \ k\ge 0.
$$
\[
\label{xFGpmDef} F(c,\t_0,\t_1,...)={1\/|\dot{W}(\s_0)|}\,,\qquad
G_\pm(c,\t_0,\t_1,\vk_1,\t_2,\vk_2,...)= -1+\sum_{n\ge
1}\lt({e^{\pm\vk_n}\/|\dot{W}(\s_n)|}-2\rt).
\]
Note that the sum in the definition of $G_\pm$ converges due to Theorem
\ref{TrubowitzNonFixAB} and Proposition \ref{PropIdentitiesAB}. Using these notations, we
rewrite the identities (\ref{IdentitiesAB}) in the following form:
\[
\label{xIdentsAB} e^{-\vk_0}F+G_-=-a,\qquad e^{\vk_0}F+G_+=-b.
\]
\begin{lemma}
\label{LmMonotMM} (i) Let $c\in\R$, $\{\t_n\}_{0}^\iy\in\cT$ and $\{\vk_n\}_{0}^\iy\in\el2_1$.
Then,
\[
\label{xFpm} {\pa F\/\pa\t_0}>0,\qquad \lim_{\t_0\to-\iy}F=0,\quad
\lim_{\t_0\to\s_1^0-\s_0^0+\t_1}F=+\iy;
\]
\[
\label{xHpm} {\pa G_\pm\/\pa\t_0}>0,\qquad \lim_{\t_0\to-\iy}G_\pm=-\iy,\quad
\lim_{\t_0\to\s_1^0-\s_0^0+\t_1}G_\pm=+\iy.
\]
(ii) For each $a,b,c\in\R$, $\{\t_{n+1}\}_{0}^\iy\in\cT^{(1)}$,
$\{\vk_{n+1}\}_{0}^\iy\in\el2_1$ there exist unique $\t_0<\s_1^0-\s_0^0+\t_1$ and $\vk_0\in\R$
such that identities (\ref{xIdentsAB}) (or, equivalently, (\ref{IdentitiesAB})) hold true.
\end{lemma}
\begin{proof} (i) The proof repeats the proof of Lemma \ref{LmMonot}.\\
(ii) Note that identities (\ref{xIdentsAB}) are equivalent to
$$
F^2=(a+G_-)(b+G_+),\qquad e^{\vk_0}=-{F\/a+G_-}=-{b+G_+\/F}.
$$
We will consider $F$ and $G_\pm$ as functions of $\t_0$. It follows from (\ref{xFpm}) that the
function $F(\t_0)$ is positive for all $\t_0<\s_1^0-\s_0^0+\t_1$. Hence, $a+G_-(\t_0)<0$ and
$b+G_+(\t_0)<0$. Due to (\ref{xHpm}), there exist unique points $\t_\pm^*<\s_1^0-\s_0^0+\t_1$
such that $a+G_-(\t_-^*)=b+G_+(\t_+^*)=0$. Since the functions $G_\pm(\t_0)$ are increasing,
we have $\t_0<\t_\pm^*$. Put $\t^*=\min(\t_-^*,\t_+^*)$ and
$$
G(\t_0)=(a+G_-(\t_0))\cdot(b+G_+(\t_0)),\quad \t_0\le\t_*.
$$
Using (\ref{xHpm}), it is easy to see that
$$
G'(\t_0)<0,\quad \lim_{\t_0\to-\iy}G(\t_0)=+\iy,\quad G(\t_*)=0;\qquad (F^2)'(\t_0)>0,\quad
\lim_{\t_0\to-\iy}(F(\t_0))^2=0.
$$
Therefore, the equation $(F(\t_0))^2=G(\t_0)$ has the unique solution $\t_0\in(-\iy,\t_*)$ and
\linebreak $\vk_0=\log|{F(\t_0)\big/(a\!+\!G_-(\t_0))}|=\log|{(b\!+\!G_+(\t_0))\big/F(\t_0)}]$
is uniquely determined by $\t_0$.
\end{proof}

\begin{proof}[{\bf Proof of Proposition \ref{CorFixedAB}.}]
Fix some $a,b\in\R$ and $m\ge 0$. Due to Theorem \ref{TrubowitzNonFixAB} and Proposition
\ref{PropIdentitiesAB}, the mapping
$
q\mapsto\P(q,a,b)=\left(Q_0+2a+2b;\{\t_n(q,a,b)\}_{0}^\iy;\{\vk_n(q,a,b)\}_{0}^\iy\right)
$
is a bijection between $L^2(0,1)$ and the set
$$
\cS^{a,b}=\lt\{(c;\{\t_n\}_0^\iy;\{\vk_n\}_0^\iy)\in\R\ts\cT\ts\el2_1:\
e^{-\vk_0}F\!+\!G_-=-a\ \ {\rm and}\ \ e^{\vk_0}F\!+\!G_+=-b\rt\},
$$
where $F$ and $G_\pm$ are given by (\ref{xFGpmDef}).

\pagebreak

\no  Define the mapping $P:\cS^{a,b}\ss\R\ts\cT\ts\el2_1\to\R\ts\cT^{(1)}\ts\el2_1$ by
$$
(c;\{\t_n\}_0^\iy;\{\vk_n\}_0^\iy)\mapsto (c;\{\t_{n+1}\}_0^\iy;\{\vk_n\}_{0}^\iy).
$$
Note that $\P_{a,b}(q)=P(\P(q,a,b))$. In particular, $\P_{a,b}$ is real-analytic as a
composition of real-analytic mappings. Due to Lemma \ref{LmMonotMM} (ii), for each
$(c;\{\t_{n+1}\}_0^\iy;\{\vk_n\}_{0}^\iy)\in\R\ts\cM^{(1)}\ts\el2_1$ there exist unique
$\t_0<\s_1^0-\s_0^0+\t_1$ and  $\vk_0\in\R$ such that both equations $e^{-\vk_0}F\!+\!G_-=-a$
and $e^{\vk_0}F\!+\!G_+=-b$ hold true. Therefore, the mapping $P$ is a bijection between
$\cS^{a,b}$ and $\R\ts\cM^{(1)}\ts\el2_1$. Thus, $\P_{a,b}$ is a bijection between $L^2(0,1)$
and $\R\ts\cM^{(1)}\ts\el2_1$.

In order to prove that $\P_{a,b}^{-1}$ is real-analytic, we deduce from Theorem
\ref{TrubowitzNonFixAB} that the mapping
$$
\left(c;\{\t_n\}_{0}^\iy;\{\vk_n\}_{0}^\iy\right)\mapsto (-a;-b)=(e^{-\vk_0}F\!+\!G_-\,;
e^{\vk_0}F\!+\!G_+)
$$
is real-analytic. Lemma \ref{LmMonotMM} (i) gives ${\pa a\/\pa\t_0}\!<\!0$, ${\pa b\/\pa
\t_0}\!<\!0$, ${\pa a\/\pa\vk_0}\!=\!e^{-\vk_0}F\!>\!0$, ${\pa
b\/\pa\vk_0}\!=\!e^{\vk_0}F\!<\!0$. Hence,
$$
\det\left(\begin{array}{cc}{\pa a\big/\pa\t_0} & {\pa b\big/\pa\t_0} \cr {\pa a\big/\pa\vk_0}
& {\pa b\big/\pa\vk_0}\end{array}\right)>\,0.
$$
Then, the Implicit Function Theorem yields that $P^{-1}:\R\ts\cT^{(1)}\ts\el2_1\to\cS^{a,b}$
is real-analytic. We obtain that $(\F^b)^{-1}$ is real-analytic as a composition of
real-analytic mappings.
\end{proof}

\end{document}